\documentclass[11pt,leqno]{amsart}
\pdfoutput=1
\usepackage{amsmath, amssymb}
\usepackage{graphicx,color,hyperref}
\usepackage{verbatim,enumitem,ifpdf}

\def\smallskip{\vskip4pt plus1pt minus1pt}
\def\medskip{\vskip7pt plus2pt minus1pt}
\def\bigskip{\vskip11pt plus2pt minus1pt}

\catcode`\@=11
\newcount\@tempcnta \newcount\@tempcntb 
\def\timeofday{{%
\@tempcnta=\time \divide\@tempcnta by 60 \@tempcntb=\@tempcnta
\multiply\@tempcntb by -60 \advance\@tempcntb by \time
\ifnum\@tempcntb > 9 \number\@tempcnta:\number\@tempcntb
\else\number\@tempcnta:0\number\@tempcntb\fi}}

\def\openup{\afterassignment\@penup\dimen@=}
\def\@penup{\advance\lineskip\dimen@
  \advance\baselineskip\dimen@
  \advance\lineskiplimit\dimen@}
\newdimen\jot \jot=3pt
\newskip\plaincentering \plaincentering=0pt plus 1000pt minus 1000pt
\def\ialign{\everycr{}\tabskip\z@skip\halign}
\def\plaineqalign#1{\null\,\vcenter{\openup\jot\m@th
      \ialign{\strut\hfil$\displaystyle{##}$&$\displaystyle{{}##}$\hfil
      \crcr#1\crcr}}\,}
\newif\ifdt@p
\def\displ@y{\global\dt@ptrue\openup\jot\m@th
  \everycr{\noalign{\ifdt@p \global\dt@pfalse \ifdim\prevdepth>-1000\p@
      \vskip-\lineskiplimit \vskip\normallineskiplimit \fi
      \else \penalty\interdisplaylinepenalty \fi}}}
\def\@lign{\tabskip\z@skip\everycr{}} 
\def\displaylines#1{\displ@y \tabskip\z@skip
  \halign{\hbox to\displaywidth{$\@lign\hfil\displaystyle##\hfil$}\crcr
    #1\crcr}}
\def\plaineqalignno#1{\displ@y \tabskip\plaincentering
  \halign to\displaywidth{\hfil$\@lign\displaystyle{##}$\tabskip\z@skip
    &$\@lign\displaystyle{{}##}$\hfil\tabskip\plaincentering
    &\llap{$\@lign##$}\tabskip\z@skip\crcr
    #1\crcr}}
\def\plainleqalignno#1{\displ@y \tabskip\plaincentering
  \halign to\displaywidth{\hfil$\@lign\displaystyle{##}$\tabskip\z@skip
    &$\@lign\displaystyle{{}##}$\hfil\tabskip\plaincentering
    &\kern-\displaywidth\rlap{$\@lign##$}\tabskip\displaywidth\crcr
    #1\crcr}}
\def\plaincases#1{\left\{\,\vcenter{\normalbaselines\m@th
    \ialign{$##\hfil$&\quad##\hfil\crcr#1\crcr}}\right.}
\def\plainmatrix#1{\null\,\vcenter{\normalbaselines\m@th
    \ialign{\hfil$##$\hfil&&\quad\hfil$##$\hfil\crcr
      \mathstrut\crcr\noalign{\kern-\baselineskip}
      #1\crcr\mathstrut\crcr\noalign{\kern-\baselineskip}}}\,}

\newdimen\plainitemindent \plainitemindent=6mm

\def\plainitem#1{\par\noindent\hangindent\plainitemindent%
            \rlap{#1}\kern\plainitemindent\ignorespaces}

\catcode`\@=12

\font\bigbf=cmbx10 at 12pt

\font\tenCal=eusm10
\font\sevenCal=eusm7
\font\fiveCal=eusm5
\newfam\Calfam
  \textfont\Calfam=\tenCal
  \scriptfont\Calfam=\sevenCal
  \scriptscriptfont\Calfam=\fiveCal
\def\Cal{\fam\Calfam\tenCal}

\font\tenBbb=msbm10
\font\sevenBbb=msbm7
\font\fiveBbb=msbm5
\newfam\Bbbfam
  \textfont\Bbbfam=\tenBbb
  \scriptfont\Bbbfam=\sevenBbb
  \scriptscriptfont\Bbbfam=\fiveBbb
\def\Bbb{\fam\Bbbfam\tenBbb}

\def\bB{{\Bbb B}}
\def\bC{{\Bbb C}}
\def\bN{{\Bbb N}}
\def\bP{{\Bbb P}}
\def\bQ{{\Bbb Q}}
\def\bR{{\Bbb R}}

\def\cO{{\Cal O}}

\def\rank{\mathop{\rm rank}}
\def\Id{{\rm Id}}
\def\tr{\mathop{\rm tr}\nolimits}
\def\Hom{\mathop{\rm Hom}\nolimits}
\def\Herm{\mathop{\rm Herm}\nolimits}
\def\MAVol{\mathop{\rm MAVol}\nolimits}
\def\det{\mathop{\rm det}\nolimits}
\def\bu{\hbox{$\scriptstyle\bullet$}}

\let\ol=\overline
\def\square{{\hfill \hbox{
\vrule height 1.453ex  width 0.093ex  depth 0ex
\vrule height 1.5ex  width 1.3ex  depth -1.407ex\kern-0.1ex
\vrule height 1.453ex  width 0.093ex  depth 0ex\kern-1.35ex
\vrule height 0.093ex  width 1.3ex  depth 0ex}}}
\def\qed{$\square$}

\ifpdf
\def\Nearrow{\mathop{\hbox{\pdfsetmatrix{0.809017 0.587785 -0.587785 0.809017}%
\rlap{\hbox{\kern5pt$\Longrightarrow$}}\pdfsetmatrix{0.809017 -0.587785
0.587785 0.809017}}\kern18pt}}
\def\Searrow{\mathop{\hbox{\pdfsetmatrix{0.809017 -0.587785 0.587785 0.809017}%
\rlap{\hbox{\kern1pt$\Longrightarrow$}}\pdfsetmatrix{0.809017 0.587785
  -0.587785 0.809017}}\kern18pt}}
\else
\def\Nearrow{\mathop{
}\kern18pt}
\def\Searrow{\mathop{
}\kern18pt}
\fi

\def\today{\ifcase\month\or
January\or February\or March\or April\or May\or June\or July\or August\or
September\or October\or November\or December\fi \space\number\day,
\number\year}

\def\Bibitem#1&#2&#3&#4&%
{\hangindent=2cm\hangafter=1
\noindent\rlap{\hbox{\bf #1}}\kern2cm{\rm #2}{\it #3}{\rm #4.}\smallskip}

\title[Monge-Amp\`ere functionals for a holomorphic vector bundle]
{Monge-Amp\`ere functionals for the curvature\vskip3pt
tensor of a holomorphic vector bundle}
\author[Jean-Pierre Demailly]{Jean-Pierre Demailly}

\begin{document}

\begin{abstract} Let $E$ be a holomorphic vector bundle on a projective
manifold $X$ such that $\det E$ is ample. We introduce three functionals
$\Phi_P$ related to Griffiths, Nakano and dual Nakano positivity
respectively. They can be used to define new concepts of volume for
the vector bundle $E$, by means of generalized Monge-Amp\`ere integrals
of~$\Phi_P(\Theta_{E,h})$, where $\Theta_{E,h}$ is the
Chern curvature tensor of $(E,h)$. These volumes are shown to satisfy
optimal Chern class inequalities. We also prove that the functionals $\Phi_P$
give rise in a natural way to elliptic differential systems of
Hermitian-Yang-Mills type for the curvature, in such a way that the
related $P$-positivity threshold of $E\otimes(\det E)^t$, where
$t>-1/\rank E$, can possibly be investigated by studying the infimum of
exponents $t$ for which the Yang-Mills differential system has a solution.


\vskip3pt\noindent
{\sc MSC Classification 2020.} 32J25, 53C07

\vskip3pt\noindent
{\sc Funding.} The author has been supported by the Advanced ERC grant
ALKAGE, no 670846 from September 2015, attributed by the
European Research Council.
\end{abstract}

\maketitle

\hbox to\textwidth{\hfill\it Dedicated to Professor L\'aszl\'o
Lempert on the occasion of his 70${}^{\,th}$ birthday}
\vskip1cm

\noindent
{\bigbf 1. Introduction}\medskip

Let $X$ be a projective $n$-dimensional manifold, and $E\to X$ a
holomorphic vector bundle equipped with a smooth hermitian metric~$h$.
Putting $\rank E=r$, the Chern curvature tensor
$\Theta_{E,h}=i\nabla_{E,h}^2$ can be written
$$
\Theta_{E,h}=i\sum_{1\le j,k\le n,\,1\le\lambda,\mu\le r}c_{jk\lambda\mu}dz_j\wedge
d\ol z_k\otimes e_\lambda^*\otimes e_\mu\leqno(1.1)
$$
in terms of holomorphic coordinates $(z_1,\ldots,z_n)$ on $X$ and of
an orthonormal frame $(e_\lambda)_{1\le\lambda\le r}$ of $E$. There is an
associated quadratic form $\widetilde\Theta_{E,h}$ on $T_X\otimes E$ defined by
$$
\widetilde\Theta_{E,h}(\gamma):=\sum_{1\le j,k\le n,\,1\le\lambda,\mu\le r}
c_{jk\lambda\mu}\gamma_{j\lambda}\ol\gamma_{k\mu},\quad
\gamma=\sum_{j,\lambda}
\gamma_{j\lambda}{\partial\over\partial z_j}\otimes e_\lambda\in T_X\otimes E,
\leqno(1.2)
$$
so that we have in particular
$$
\widetilde\Theta_{E,h}(\xi\otimes v):=
\langle\Theta_{E,h}(\xi,\ol\xi)\cdot v,v\rangle_h=
\sum_{1\le j,k\le n,\,1\le\lambda,\mu\le r}
c_{jk\lambda\mu}\xi_j\ol\xi_k v_\lambda\ol v_\mu.
$$
As is well known, the dual hermitian bundle $(E^*,h^*)$ has a curvature tensor
that is the opposite of the transpose of $\Theta_{E,h}$, and
for $\gamma=\sum_{j,\lambda}\gamma_{j\lambda}
{\partial\over\partial z_j}\otimes e^*_\lambda\in T_X\otimes E^*$ we have
$$
\plainleqalignno{
-\Theta_{E^*,h^*}&={}^T\Theta_{E,h}=
\sum_{1\le j,k\le n,\,1\le\lambda,\mu\le r}c_{jk\mu\lambda}dz_j\wedge
d\ol z_k\otimes (e^*_\lambda)^*\otimes e^*_\mu,&(1.3)\cr
-\widetilde\Theta_{E^*,h^*}(\gamma)&={}^T\widetilde\Theta_{E,h}(\gamma)=
\sum_{1\le j,k\le n,\,1\le\lambda,\mu\le r}
c_{jk\mu\lambda}\gamma_{j\lambda}\ol\gamma_{k\mu}>0.
&(1.4)\cr}
$$
Let us recall the following standard positivity concepts.
\medskip

\noindent
{\bf 1.5. Definition.} {\it The hermitian bundle $(E,h)$ is said to be
{\plainitemindent=6.5mm\par
\plainitem{\rm(a)} Griffiths positive if $\,\widetilde\Theta_{E,h}(\xi\otimes v)>0$
for all decomposable nonzero tensors $\xi\otimes v\in T_X\otimes E,$
\plainitem{\rm(b)} Nakano positive if $\,\widetilde\Theta_{E,h}
(\gamma)>0$ for all nonzero tensors $\gamma\in T_X\otimes E,$
\plainitem{\rm(c)} dual Nakano positive if $\,{}^T\widetilde\Theta_{E,h}(\gamma)>0$
for all nonzero tensors $\gamma\in T_X\otimes E^*$.
\vskip2pt}
\noindent
One says that $E$ itself has one of these three positivity properties if
it possesses a smooth hermitian metric $h$ satisfying the corresponding
positivity assumption.}
\medskip

\noindent This definition gives rise to well known implications
$$
\plainleqalignno{
E ~\hbox{Nakano positive}\,\Searrow&\cr
&~E~\hbox{Griffiths positive}~{}\Longrightarrow E~\hbox{ample}.\cr
E~\hbox{dual Nakano positive}\,\Nearrow&\cr}
$$
The last implication comes from the Kodaira embedding theorem [Kod54] and
the easy verification that the Griffiths positivity of $\Theta_{E,h}$
implies the positivity of the curvature of the
induced metric on the tautological line bundle $\cO_{\bP(E)}(1)$, where
$\bP(E)$ is the projectivized bundle of hyperplanes of~$E$. A basic
problem raised by [Gri69, Problem (0.9)] is

\medskip\noindent
{\bf 1.6. Griffiths problem.} {\it Does it hold that~
$E~\hbox{ample}\Rightarrow E~\hbox{Griffiths positive}~?$}
\medskip

\noindent One might wonder whether the ampleness of $E$ would
even imply the Nakano or dual Nakano positivity of $E$, but it turns
out that none of these implications holds true. In fact the tangent
bundle $E=T_X$ of the complex projective space $X=\bP^n$ is ample but
not Nakano positive (the fact that
$H^{n-1,n-1}(\bP^n,\bC)=H^{n-1}(X,K_X\otimes T_X)\neq 0$ would
contradict the Nakano vanishing theorem [Nak55]), and the cotangent
bundle $E=T^*_X$ of a ball quotient $X=\bB^n/\Gamma$ is ample but not
dual Nakano positive, since $\Id_E\in H^0(X,\Omega^1_X\otimes E^*)\neq
0$ contradicts the dual version of the Nakano vanishing theorem
([Nak55], [DemB]). The latter fact, that was very briefly mentioned in
[LSY13, p.~304], had been overlooked in [Dem21], where we proposed an
approach to investigate the dual Nakano positivity of an ample vector
bundle. However, as we will see here, the above counterexample
raises new interesting problems and does not invalidate the approach of [Dem21].
We thank Dr Junsheng Zhang for pointing out to us the observation
made in [LSY13].

In section~2, we introduce three functionals $\Phi_N(\Theta_{E,h})$,
$\Phi_{N^*}(\Theta_{E,h})$, $\Phi_G(\Theta_{E,h})$ and corres\-ponding
integrated Monge-Amp\`ere volumes $\MAVol_N(E)$, $\MAVol_{N^*}(E)$,
$\MAVol_G(E)$ that are related respectively to Nakano, dual Nakano and Griffiths
positivity. One can check that these Monge-Amp\`ere volumes
reach their maximum value if and only if the bundle $E$ is projectively
flat -- see Corollary~2.7 for a detailed statement. The corresponding
densities are determinants of the curvature tensor that can be used to
define global scalar equations for the curvature.
In section 3, extending the approach proposed in [Dem21], we
show that it suffices to add a trace free Hermite-Einstein condition
to a scalar determinantal equation to yield families of elliptic systems
of Yang-Mills type, denoted respectively ${\rm YM}_{N,\beta}(t)$,
${\rm YM}_{N^*,\beta}(t)$, ${\rm YM}_{G,\beta}(t)$, depending on a
time parameter~$t$ and a suitable positive constant $\beta$,
where the unknown is a time dependent hermitian metric $h_t$ on $E$.
These solutions could hopefully help in the study of the Griffiths' problem
if one could obtain an appropriate existence theorem.
On a more differential geometric side, if $(E,h)$ is a hermitian
vector bundle and $t\in\bR$ a real number, we consider formally the
curvature tensor of $E\otimes (\det E)^t$, namely
$$
\Theta_{E,h}+t\,\Theta_{\det E,\det h}\otimes\Id_E=
\Theta_{E,h}+t\,\tr_E\Theta_{E,h}\otimes\Id_E,
\leqno(1.7)
$$
where $\Theta_{\det E,\det h}$ is the $(1,1)$-curvature form of the
determinant bundle $\det E=\Lambda^rE$ and $\tr_E$ the trace
operator on $\Hom(E,E)$. We introduce the following
threshold values, defined for any vector bundle possessing an ample
determinant.
\medskip

\noindent
{\bf 1.8. Definition.} {\it Let $E\to X$ be a holomorphic vector bundle
such that $\det E$ is ample. We define the Nakano, dual Nakano, Griffiths
and ample thresholds, denoted respectively
$$
\tau_N(E),\quad \tau_{N^*}(E),\quad\tau_G(E),\quad\tau_A(E),
$$
to be the infimum of values $t\in\bR$ such that there exists a smooth
hermitian metric $h$ for which
$\Theta_{E,h}+t\,\Theta_{\det E,\det h}\otimes\Id_E$ is
Nakano, dual Nakano, Griffiths positive, respectively the infimum of
$t\in\bQ$ such that $E\otimes (\det E)^t$ is ample $($i.e.\
$\cO_{\bP(E)}(1)\otimes\pi^*(\det E)^{t}$ is a $\bQ$-ample line bundle
on the total space of the projectivized bundle $\pi:\bP(E)\to X)$.}
\medskip

\noindent Notice that Nakano and dual Nakano positivity are
stronger than Griffiths positivity, the latter being itself stronger than
ampleness, hence we always have
$$
\tau_N(E)\geq\tau_G(E)\geq\tau_A(E),\quad \tau_{N^*}(E)\geq\tau_G(E)
\geq\tau_A(E).\leqno(1.9)
$$
Since $E\otimes(\det E)^{-1/r}$ has trivial determinant, no positivity
property can hold for it, and we conclude that $\tau_A(E)\geq -1/r$.
The equality may however occur, e.g.\ when $E=L^{\oplus r}$ is
the direct sum of $r$ copies of an ample line bundle~$L$. In fact,
the following fact is easy to prove.
\medskip

\noindent
{\bf 1.10. Proposition.} {\it Let $E\to X$ be a holomorphic vector bundle
such that $\det E$ is ample. If~$\tau_G(E)=-1/r$, then
$F=E\otimes(\det E)^{-1/r}$ is numerically flat, in other words, as
a $\bQ$-vector bundle, we have $E=F\otimes L$ where $F$ is numerically
flat of rank $r$, and $L=(\det E)^{1/r}$ is an ample $\bQ$-line bundle
$($one will refer to this situation by saying
that $E$ is projectively numerically flat, see Definition~$2.8)$.
Then~we have
$$
\tau_N(E)=\tau_{N^*}(E)=\tau_G(E)=\tau_A(E)=-{1\over r}.
$$}
\vskip-2pt

\noindent
In this setting, the Griffiths problem translates into the conjectural
implication
$$
\hbox{$E$ ample $\Rightarrow\tau_G(E)<0$~?}\leqno(1.11)
$$
On the other hand, in the counterexamples $E=T_{\bP^n}$
(resp.\ $E=T^*_{\bB^n/\Gamma}$) just mentioned for Nakano
(resp.\ dual Nakano) positivity, one can check that
$E$ is in fact Griffiths positive and Nakano (resp.\ dual Nakano)
semipositive, hence we have $\tau_G(E)<0$, while
$\tau_N(E)=0$ (resp.~$\tau_{N^*}(E)=0$). By investigating the curvature
of direct images of adjoint positive line bundles, Berndtsson [Ber09]
(see also Mourougane-Takayama [MoT07]) has proved that
$$
\hbox{
$E$ ample $\Rightarrow S^mE\otimes\det E=\pi_*(\cO_{\bP(E)}(m)\otimes K_{\bP(E)})$
Nakano positive for all $m\in\bN$.}
\leqno(1.12)
$$
We infer from this that $E$ ample implies $\tau_N(E)<1$
and $\tau_N(S^mE)<{r\over r_m}$ where $r_m$ is the rank
of $S^mE$, namely $r_m={m+r-1\choose r-1}$. Furthermore, we know
that $S^mE$ generates its jets for $m\geq m_0$ large, hence
$\tau_{N^*}(S^mE)<0$ for $m\geq m_0$. In [LSY13, Cor.~4.12], it is further
proved that
$$
\hbox{
$E$ ample $\Rightarrow S^mE\otimes\det E$ dual Nakano positive
for all $m\in\bN$,}
\leqno(1.13)
$$
hence we have as well $\tau_{N^*}(E)<1$ and $\tau_{N^*}(S^mE)<{r\over r_m}$
when $E$ is ample. The only counter-examples we know about still leave room
for the following question.

\medskip\noindent
{\bf 1.14. Question.} {\it Assume that $E$ is an ample vector bundle.
Are there examples for which $\tau_N(E)>0$, $\tau_{N^*}(E)>0$ or
$\tau_G(E)\geq 0~?$}
\medskip\noindent
Of course finding  an example with $\tau_G(E)\geq 0$ would be equivalent
to answer negatively Griffiths' problem 1.6. On the PDE side, 
our main result is as follows (see section~\S$\,$3).
\medskip

\noindent
{\bf 1.15. Theorem.} {\it Given any value $t_0$ such that
$E\otimes (\det E)^{t_0}>_P0$, one can always arrange the
corresponding differential systems
${\rm YM}_{N,\beta}(t)$, ${\rm YM}_{N^*,\beta}(t)$, ${\rm YM}_{G,\beta}(t)$ to
be elliptic inver\-tible and to have unique solutions that depend
continuously $($and even differentiably$)$ on $t$ on a small interval
$[t_0-\delta_0,t_0+\delta_0]$, $\delta_0>0$.}
\medskip

\noindent
The proof depends only on the theory of elliptic equations and on
the implicit function  theorem. In the end, checking the
ellipticity is just a sophisticated exercise of linear algebra.
A natural problem is whether such Yang-Mills type equations can be
used to compute the positivity thresholds, by trying to get solutions
for $t\in{}]-1/r,t_0]$ as small as possible.
\medskip

\noindent
{\bf 1.16. Question.} {\it Can one design the Yang-Mills systems
${\rm YM}_{N,\beta}(t)$,  ${\rm YM}_{N^*,\beta}(t)$, ${\rm YM}_{G,\beta}(t)$
so that the infimum of times $t_{\inf}$ for which
a smooth solution exists on $]t_{\inf},t_0]$
coincides respectively with the positivity thresholds
$\tau_N(E)$,~$\tau_{N^*}(E)$, $\tau_G(E)$, for suitably
chosen initial data at $t=t_0$ $($or whatever they are$\,)~?$}
\medskip

\noindent
Getting $t_{\inf}=t_{\inf}(\beta)$ to converge to the positivity threshold
$\tau_P(E)$ as $\beta\to+\infty$ instead of being equal to $\tau_P(E)$
would be good as well.
In the above question, we somehow expect that the differential
systems can be made invertible elliptic throughout an almost
maximal interval $[t_P(E)+\delta,t_0]$, $0<\delta\ll 1$, and not just
on a small interval $[t_0-\delta_0,t_0]$.
\medskip
\noindent
{\bf Acknowledgment.} We warmly thank Vamsi Pritham Pingali and Junsheng Zhang
for enlightening exchanges that contributed to clarify some of our ideas.
The reader is referred to [Pin20] and [Pin21] for related results.
\bigskip

\noindent
{\bigbf 2. Monge-Amp\`ere functionals for vector bundles}\medskip

Let $E\to X$ be a holomorphic vector bundle equipped with a smooth hermitian
metric~$h$. If the Chern curvature tensor $\Theta_{E,h}$ is Nakano
positive, then the ${1\over r}$-power of
the $(n\times r)$-dimensional determinant of the corres\-ponding
hermitian quadratic form on $T_X\otimes E$ can be seen as a
positive $(n,n)$-form
$$
\Phi_N(\Theta_{E,h})=\det_{T_X\otimes E}(\Theta_{E,h})^{1/r}
=\det(c_{jk\lambda\mu})_{(j,\lambda),(k,\mu)}^{1/r}\,
i dz_1\wedge d\ol z_1\wedge\ldots\wedge i dz_n\wedge d\ol z_n.
\leqno(2.1)
$$
Moreover, this $(n,n)$-form does not depend on the choice of coordinates
$(z_j)$ on~$X$, nor on the choice of the orthonormal frame $(e_\lambda)$ on $E$
(but $(e_\lambda)$ must be orthornormal). Similarly,
if the Chern curvature tensor $\Theta_{E,h}$ is dual
Nakano positive, we can consider the $(n\times r)$-dimensional
determinant of the hermitian quadratic form on $T_X\otimes E^*$, namely
$$
~~\Phi_{N^*}(\Theta_{E,h})=\det_{T_X\otimes E^*}({\,}^T\Theta_{E,h})^{1/r}
=\det(c_{jk\mu\lambda})_{(j,\lambda),(k,\mu)}^{1/r}\,
i dz_1\wedge d\ol z_1\wedge\ldots\wedge i dz_n\wedge d\ol z_n,
\leqno(2.2)
$$
and view it as a positive $(n,n)$-form. Finally, if $\Theta_{E,h}$ is
Griffiths positive, the most natural substitute for (2.1) and (2.2) is
$$
~~\Phi_G(\Theta_{E,h})
=\inf_{\strut{\scriptstyle v\in E},\atop\scriptstyle|v|_h=1}\;{1\over n!}\,
\big(\langle\Theta_{E,h}\cdot v,v\rangle_h\big)^n
=\inf_{|v|=1}\;{1\over n!}\,\Bigg(\sum_{1\leq \lambda,\mu\leq r}
c_{jk\lambda\mu}\,v_\lambda\overline v_\mu\,
i dz_j\wedge d\ol z_k\Bigg)^{\!n}.\leqno(2.3)
$$
It is easy to see that the three volume forms coincide when $(E,h)$ is
projectively flat, namely when $\Theta_{E,h}=\alpha\otimes\Id_E$ where
$\alpha$ is a positive $(1,1)$-form on $X$ (which is then
equal to ${1\over r}\tr_E\Theta_{E,h}={1\over r}\Theta_{\det E,\det h}$).
In this case, we clearly have
$$
\Phi_N(\Theta_{E,h})=\Phi_{N^*}(\Theta_{E,h})
=\Phi_G(\Theta_{E,h})={1\over n!}\,\alpha^n
={1\over n!\, r^n}\,(\Theta_{\det E,\det h})^n.\leqno(2.4)
$$
In general, we have the following inequalities.
\medskip

\noindent
{\bf 2.5. Proposition.} {\it Let $(E,h)$ be a hermitian vector bundle.
{\plainitemindent=6.5mm\vskip2pt
\plainitem{\rm(a)} If $\Theta_{E,h}$ is Nakano positive, then
$\displaystyle\Phi_N(\Theta_{E,h})\leq
{1\over n!\, r^n}\,(\Theta_{\det E,\det h})^n$.
\plainitem{\rm(b)} If $\Theta_{E,h}$ is dual Nakano positive, then
$\displaystyle\Phi_{N^*}(\Theta_{E,h})\leq
{1\over n!\, r^n}\,(\Theta_{\det E,\det h})^n$.

\plainitem{\rm(c)} If $\Theta_{E,h}$ is Griffiths positive, then
$\displaystyle\Phi_G(\Theta_{E,h})\leq
{1\over n!\, r^n}\,(\Theta_{\det E,\det h})^n$.
\vskip6pt\noindent
In all three cases, the equality of volume forms occurs if and only if
$(E,h)$ is projectively flat and~$\Theta_{\det E,\det h}>0$.}}
\medskip

\noindent{\it Proof.} (a) We take $h$ to be a hermitian metric on $E$ such that
$\Theta_{E,h}$ is Nakano positive, and consider the K\"ahler metric
$$
\omega=\Theta_{\det E,\det h}=\tr_E\Theta_{E,h}.
$$
If $(\alpha_j)_{1\le j\le nr}$ are the eigenvalues of the associated
hermitian form $\widetilde\Theta_{E,h}$ with respect to
$\omega\otimes h$, we have
$$
\det_{T_X\otimes E^*}\big(\Theta_{E,h}\big)^{1/r}
=\bigg(\prod_j\alpha_j\bigg)^{1/r}{\omega^n\over n!}
$$
and $\big(\prod_j\alpha_j\big)^{1/nr}\le
{1\over nr}\sum_j\alpha_j$ by the inequality between the geometric and
arithmetic means. Since
$$
\sum_j\alpha_j=\tr_\omega\big(\tr_E\Theta_{E,h}\big)=\tr_\omega\omega=n,
$$
we obtained the asserted inequality
$$
\det_{T_X\otimes E}\big(\Theta_{E,h}\big)^{1/r}
\leq \bigg({1\over nr}\sum_j\alpha_j\bigg)^n{\omega^n\over n!}
={1\over n!\,r^n}\,\omega^n.
$$
(b) In the case of dual Nakano positivity, the proof is almost identical,
except that we take the $\alpha_j$'s to be the eigenvalues of
${}^T\widetilde\Theta_{E,h}$ with respect to $\omega\otimes h^*$ on
$T_X\otimes E^*$. In both cases, the equality of volume forms occurs
if only if all eigenvalues $\alpha_j$ are equal at all points, and then
we must have $\alpha_j={1\over r}$, hence
$\Theta_{E,h}={1\over r}\,\omega\otimes h$.
\medskip

\noindent
(c) When $(E,h)$ is Griffiths positive, we pick an $h$-orthonormal frame
$(e_\lambda)_{1\leq\lambda\leq r}$ of $E$ and observe that we have by definition
$$
\Phi_G(\Theta_{E,h})\leq {1\over n!}\,\big(\langle\Theta_{E,h}\cdot
e_\lambda,e_\lambda\rangle_h\big)^n,\quad 1\leq\lambda\leq r.\leqno(*)
$$
We view this inequality as a comparison between positive real numbers by
referring to the volume form $dV={1\over n!}\omega^n$ of the
metric $\omega=\Theta_{\det E,\det h}$.
Let us consider the $(1,1)$-form
$A_\lambda=\langle\Theta_{E,h}\cdot e_\lambda,e_\lambda\rangle_h$
on $(T_X,\omega)\simeq (\bC^n,{\rm std})$ as a hermitian form (or matrix) on
$\bC^n$. Then we have ${1\over n!}\,
\langle\Theta_{E,h}\cdot e_\lambda,e_\lambda\rangle_h^n=\det A_\lambda\,dV$.
We use the well-known fact that the
function $A\mapsto(\det A)^{1/n}$ is concave on the cone of
positive hermitian $(n\times n)$-matrices. This implies
$$
{1\over r}\sum_{\lambda=1}^r(\det A_\lambda)^{1/n}\leq
\Bigg(\det\Bigg({1\over r}\sum_{\lambda=1}^r A_\lambda\Bigg)\Bigg)^{1/n},
$$
in other words
$$
{1\over r}\sum_{\lambda=1}^r\bigg(
{1\over n!}\,\langle\Theta_{E,h}\cdot e_\lambda,e_\lambda\rangle_h^n\bigg)^{\!1/n}
\leq\Bigg({1\over n!}\,\Bigg({1\over r}\sum_{\lambda=1}^r
\langle\Theta_{E,h}\cdot e_\lambda,e_\lambda\rangle_h\Bigg)^{\!n}\;\Bigg)^{\!1/n}
=\bigg({1\over n!}\Big({1\over r}\omega\Big)^{\!n}\;\bigg)^{\!1/n}.
$$
By $(*)$, the left hand side is greater or equal to $\Phi_G(\Theta_{E,h})^{1/n}$,
and by taking the $n$-th power of the above inequality, we find
$\Phi_G(\Theta_{E,h})\leq {1\over n!\,r^n}\,\omega^n$,
as desired. The only line segments that lie in the graph of
$A\mapsto(\det A)^{1/n}$ project into rays of the cone of positive
hermitian matrices. Therefore, the equality case may occur only
when the hermitian forms $A_\lambda$ are proportional and we have
$\Phi_G(\Theta_{E,h})={1\over n!}\,\langle\Theta_{E,h}\cdot
e_\lambda,e_\lambda\rangle^n$ for each $\lambda$. This forces the
$(1,1)$-forms $A_\lambda$ to be equal, and therefore equal to
${1\over r}\omega$, for any choice of $h$-orthonormal frame $(e_\lambda)$.
It follows that $(E,h)$ must be projectively flat.\hfill\qed\medskip

\noindent
By considering their integrals over $X$, the above functionals give rise to
interesting concepts of volume for vector bundles.
\medskip

\noindent
{\bf 2.6. Definition.} {\it Let $E\to X$ be a holomorphic vector bundle.
If $E$ is $P$-positive, where $P$ is any of the symbols
$N,N^*$ or $G$, we define the related Monge-Amp\`ere volume of $E$ to be
$$
\MAVol_P(E)=\sup_h {1\over (2\pi)^n}\int_X\Phi_P(\Theta_{E,h}),
$$
where the supremum is taken over all smooth metrics $h$ on $E$ such that
$\Theta_{E,h}$ is $P$-positive.}
\medskip

\noindent
By Proposition~2.5, the supremum is always finite, and 
in fact we immediately get the following upper bound from the fact
that ${1\over 2\pi}\Theta_{\det E,\det h}$ is a $(1,1)$-form representing
the first Chern class $c_1(\det E)=c_1(E)$.\medskip

\noindent
{\bf 2.7. Corollary.} {\it For any $P$-positive vector bundle $E$,
we have
$$
\MAVol_P(E)\leq {1\over n!\,r^n}\,c_1(E)^n.
$$
Moreover, the equality occurs, with the supremum being a maximum, if and
only if $E$ is projectively flat.}
\medskip

\noindent
It may happen that the equality occurs for the supremum, without $E$ being
projectively flat. In fact, one has to take account the following more
general situation.
\medskip

\noindent
{\bf 2.8. Definition.} {\it We say that a rank $r$ vector bundle $E$ is
numerically projectively flat if \hbox{$F=S^rE\otimes(\det E)^{-1}$}
is numerically
flat, i.e;\ both $F$ and $F^*$ are nef vector bundles. An equivalent condition
is that the $\bQ$-line bundles $\cO_{\bP(E)}(1)\otimes\pi^*(\det E)^{-1/r}$ over
$\bP(E)$ and $\cO_{\bP(E^*)}(1)\otimes\pi^*(\det E^*)^{-1/r}$ over $\bP(E^*)$
are both nef.}
\medskip

\noindent If $\det E$ happens to admit an $r$-th root $(\det E)^{-1/r}$ that
is a genuine line bundle on $X$, then the numeric projective flatness of
$E$ is equivalent to $F=E\otimes(\det E)^{-1/r}$ being numerically flat.
In that case, we know by [DPS94, Theorem 1.18] (assuming $X$ to be projective),
that this is equivalent to the existence of a filtration
$0=F_0\subset F_1\subset\cdots\subset F_k=F$ by vector bundles $F_j$ such that
the graded pieces $F_j/F_{j-1}$ are hermitian flat for $1\leq j\leq k$, i.e.\
given by representatons of $\pi_1(X)$ into the a unitary group $U(r_j)$.
Since there always exists a finite morphism $Y\to X$ such that the pull-back of
$\det E$ to $Y$ admits an $r$-th root on $Y$, we can always obtain such a
filtration by pulling back $E$ itself. The above considerations lead
to the following fact.
\medskip

\noindent
{\bf 2.9. Proposition.} {\it Assume that $E$ is numerically projectively
flat, that $\det E$ is ample and admits an $r$-th root on $X$. Then
$$
\MAVol_N(E)=\MAVol_{N^*}(E)=\MAVol_G(E)={1\over n!\,r^n}\,c_1(E)^n.
$$}
\vskip-4pt

\noindent{\it Proof}. The proof proceeds by showing that there exist
smooth metrics $\tilde h_\varepsilon$ on the numerically flat bundle
$F=E\otimes (\det E)^{-1/r}$, such that the curvature tensor
$\Theta_{F,\tilde h_\varepsilon}$ is arbitrarily small in $L^\infty$ norm.
This is a standard fact, resulting from the property that the filtered
bundle $F$ deforms to its graded bundle
$G=\smash{\bigoplus_{1\leq j\leq k}F_j/F_{j-1}}$. In fact, it is enough to
fix a $C^\infty$ splitting of the filtration $(F_j)$ and hermitian flat
metrics $h_j$ on the graded pieces $F_j/F_{j-1}$. We can then use
$\smash{\tilde h_\varepsilon=\bigoplus_{1\leq j\leq k}\varepsilon^{k-j}h_j}$
on $F\simeq G$ (as~a $C^\infty$ vector bundle). An easy
check shows that the second fundamental forms of the Chern connections
become arbitrarily small in $L^\infty$ norms.
We take a metric of positive curvature $\eta$ on $\det E$, and consider
the metrics $h_\varepsilon=\tilde h_\varepsilon\otimes\eta^{1/r}$ on $E$.
One can then see that the supremum of the Monge-Amp\`ere integrals
over the family $(h_\varepsilon)$ reaches
the equality in~2.9. When the filtration is non split, the supremum is
never a maximum, as this would imply $E$ to be projectively flat
by Proposition~2.5.\hfill\qed
\medskip

\noindent
{\bf 2.10. Complements.} (a) The argument used in the proof of Proposition~2.9
also implies Proposition 1.10, even without assuming that
the $n$-th root $(\det E)^{1/r}$ exists on $X$. In fact, we can
extract the $n$-th root of $\det E$ by pulling back via a finite morphism
$\mu:Y\to X$. We then get a family of metrics $\tilde h$ on
$\mu^*E$ achieving the desired threshold $-1/r$ over $Y$. We define a
metric$h$  on $E$ (or rather $h^*$ on $E^*$) by putting
$$
|\xi|_{h^*(x)}^2=\sum_{y\in\mu^{-1}(x)}|\xi|_{\tilde h^*(y)}^2,\quad
\xi\in E^*_x,~~x\in X,~~y\in Y,
$$
where the sum is counted with multiplicity at branched points.
Since Griffiths sempositivity is equivalent to the plurisubharmonicity of
$|\xi|_{h^*}^2$ on $E^*$, this process preserves Griffiths (semi)-positivity;
in general the metric $h^*$ is just continuous, but we can apply a Richberg
regularization process to make it smooth. The armument is complete for the
Griffiths threshold $\tau_G$. For the Nakano and dual Nakano positivity,
we use the fact that $E$ is a subbundle of $\mu_*\mu^*E$ (and
likewise for $E^*$), Nakano seminegativity being preserved by going
to subbundles.\smallskip

\noindent
(b) In the case of a completely split bundle
$E=\bigoplus_{j=1}^rE_j$ with ample factors $E_j$ of rank~$1$,
equipped with a split metric $\smash{h=\bigoplus_{j=1}^rh_j}$,
Yau's theorem [Yau78] allows us to normalize the metrics $h_j$
to have proportional volume forms
$({1\over 2\pi}\Theta_{E_j,h_j})^n=\beta_j\omega^n$ for any K\"ahler metric
$\omega\in c_1(E)$, $\beta_j>0$ being a suitable constant.
We then get $\beta_j=c_1(E_j)^n/c_1(E)^n$, and find
$$
\plaineqalign{
{1\over (2\pi)^n}\int_X\Phi_N(\Theta_{E,h})
={1\over (2\pi)^n}\int_X\Phi_{N^*}(\Theta_{E,h})
&=\Bigg(\prod_{j=1}^r\beta_j\Bigg)^{\!1/r}\int_X{\omega^n\over n!}
={1\over n!}\,\Bigg(\prod_{j=1}^rc_1(E_j)^n\Bigg)^{\!1/r}.\cr}
$$
For $P=N,N^*$, the inequality of Corollary 2.7 then reads
$$
\Bigg(\prod_{j=1}^rc_1(E_j)^n\Bigg)^{\!{1\over r}}
\leq{1\over r^n}\,c_1(E)^n.
$$
It is an equality when $E_1=\cdots=E_r$, thus Corollary 2.7 is
optimal as far as the constant ${1\over n!\,r^n}$ is concerned.
For a completely split bundle $E=\bigoplus_{1\le j\le r}E_j$
with arbitrary ample factors, it~seems natural to conjecture
that
$$
\MAVol_N(E)=\MAVol_{N^*}(E)=
{1\over n!}\,\Bigg(\prod_{j=1}^rc_1(E_j)^n\Bigg)^{\!{1\over r}},
$$
i.e.\ that the supremum is reached for split metrics $h=\bigoplus h_j$.
In the case of the Griffiths functional, it is easy to see that
$$
\MAVol_G(E)={1\over (2\pi)^n}\int_X\Phi_G(\Theta_{E,h})
=\min_{1\leq j\leq r}\beta_j\int_X{\omega^n\over n!}
={1\over n!}\,\min_{1\leq j\leq r}c_1(E_j)^n.\leqno(2.11)
$$
In fact, $\Phi_G(\Theta_{E,h})$ is obtained by picking vectors
$v$ in the component $E_j$ for which $\beta_j$ is minimum. Moreover, for
any $G$-positive metric $h$ on $E$, even a non split one, (b1) is
proved by arguing with the induced metric $h_{|E_j}$ on $E_j$, which is again
$G$-positive as a quotient of the metric of $E$ by the projection $E\to E_j$.
\medskip

\noindent
(c) It would be interesting to characterize the ``extremal
metrics'' $h$ achieving the supremum in $\MAVol_N(E)$, $\MAVol_{N^*}(E)$,
when a maximum exists (we have seen in the proof of Proposition 2.9 that
this is not always the case). Suitable calculations (see \S3 for this)
would show that they satisfy a certain Euler-Lagrange equation
$$
\int_X
(\det \theta)^{1/r}\cdot
\tr_{T_X\otimes E^*}\Big(\theta^{-1}\cdot{\,}^T\big(
i\partial_{h^*\otimes h}\ol\partial u\big)\Big)=0\qquad
\forall u\in C^\infty(X,\Herm(E)),
\leqno(2.12)
$$
where $\theta$ is the $(n\times r)$-matrix representing
${\,}^T\Theta_{E,h}$. After integrating by parts twice, freeing
$u$ from any differentiation, we get a fourth order nonlinear
differential system that $h$ has to satisfy. Such a system is somewhat
akin to the equation for cscK metrics, in the special case~$E=T_X$.
\medskip
\noindent
(d) When $r=1$, we clearly have
$$
\Phi_N(\Theta_{E,h})=\Phi_{N^*}(\Theta_{E,h})=\Phi_G(\Theta_{E,h})=
{1\over n!}\,(\Theta_{E,h})^n,\leqno(2.13)
$$
and we infer that the integrals ${1\over(2\pi)^n}\int_X\Phi_P(\Theta_{E,h})=
{1\over n!}\,c_1(E)^n$ do not depend on $h$. In the case of ranks
$r>1$, it is natural to ask what is the infimum
$$
\inf_h{1\over (2\pi)^n}\int_X\Phi_P(\Theta_{E,h})
$$
for the three types of functionals. Let us consider the
split case $(E,h)=\bigoplus(E_j,h_j)$. By [Yau78] again, we can renormalize
$\Theta_{E_j,h_j}$ to get volume form
equalities $(2\pi)^{-n}\Theta_{E_j,h_j}^n=f_j\omega^n$ with arbitrary
functions $f_j>0$ such that $\int_X f_j\omega^n=c_1(E_j)^n$. Then
$$
\plainleqalignno{
&{1\over(2\pi)^n} \int_X\Phi_N(\Theta_{E,h})
={1\over(2\pi)^n} \int_X\Phi_{N^*}(\Theta_{E,h})
=\int_X(f_1\cdots f_r)^{1/r}\,{\omega^n\over n!},&(2.14)\cr
&{1\over(2\pi)^n}\int_X\Phi_G(\Theta_{E,h})
=\int_X\min_{1\leq j\leq r}f_j\,{\omega^n\over n!},&(2.15)\cr}
$$
and these integrals become arbitrarily small if we take the
$f_j$'s to be large on disjoint open sets, and very small elsewhere.
This example leads us to suspect that for $r>1$ and any $P=N,N^*,G$,
one always have
$$
\inf_h{1\over(2\pi)^n}\int_X\Phi_P(\Theta_{E,h})=0.\leqno(2.16)
$$

\noindent
(e) In example (2.15), we have $\Phi_G(\Theta_{E,h})=
\min_{1\leq j\leq r}f_j$. This shows that the functional $\Phi_G$ fails
to be differentiable in general, even though it is clearly
Lipschitz continuous. In order to mitigate this difficulty, it
suffices to take a large parameter $s>0$ and to consider the family
of functionals
$$
\Phi_{G,s}(\Theta_{E,h})=\Bigg(\int_{v\in E,\,|v|_h=1}
{1\over\big((\langle\Theta_{E,h}\cdot v,v\rangle_h)^n\big)^{s\phantom{!\!}}}
\;d\sigma(v)\Bigg)^{\smash{\!-1/s}}\leqno(2.17)
$$
where $d\sigma$ is the unitary invariant probability measure on the unit
sphere. Then $\Phi_{G,s}\geq\Phi_G$ and $\Phi_{G,\infty}=\lim_{s\to+\infty}
\Phi_{G,s}=\Phi_G$. A~differentiation under the integral sign shows that
$\Phi_{G,s}$ is a differentiable functional whenever $s<\infty$.
\medskip

\noindent
(f) One desirable property for our functionals $\Phi_P(\Theta_{E,h})$ is
that the conditions $\Theta_{E,h}\geq_P 0$ and $\Phi_P(\Theta_{E,h})>0$
should enforce $\Theta_{E,h}>_P 0$, thus preventing the
$P$-positivity of $\Theta_{E,h}$ to degenerate. This is clearly the
case for $\Phi_N$, $\Phi_{N^*}$, $\Phi_G$. This will be also the
case for~$\Phi_{G,s}$ for $s\geq r-1\,$: in fact, 
if $\langle\Theta_{E,h}\cdot v,v\rangle_h\ge 0$ and
$(\langle\Theta_{E,h}\cdot v,v\rangle_h)^n$ vanishes at some point
$(x_0,v_0)\in E$, $|v_0|=1$, then the differentiability of the non
negative polynomial expressing the volume form along the
$(r-1)$-dimen\-sional projectivized fiber $P(E_{x_0})$ implies
$(\langle\Theta_{E,h}\cdot v,v\rangle_h)^n/|v|^{2n}=O(|v-v_0|^2)$
near~$[v_0]\in P(E_{x_0})\,$; we conclude that the integral (2.17) is 
divergent on $P(E_{x_0})$ as well as on the unit sphere of~$E_{x_0}$,
and this shows that $\Phi_{G,s}(\Theta_{E,h})(x_0)=0$ for~$s\geq r-1$.
\bigskip

\noindent
{\bigbf 3. Hermitian-Yang-Mills equations and positivity thresholds}\medskip

Following the strategy suggested in [Dem21], we propose here to
study certain differential systems of Yang-Mills type, that
could be useful to obtain information on the positivity thresholds
of a holomorphic vector bundle. Throughout this section, we assume that
$X$ is a complex projective manifold of dimension $n$, and that $E\to X$
is a rank $r$ holomorphic vector bundle such that $\det E$ is ample.
Then there exists $t_0>0$ such that $E\otimes(\det E)^{t_0}$ has all positivity
properties $P=N,N^*,G$ we may desire. If $E$ itself is assumed to
be ample, we know by [Ber09] and [LSY13] that one can take $t_0=1$.
We consider time dependent smooth metrics $(h_t)_{t\in{}]t_1,t_0]}$
on $E$, such that
$$
\Theta_{E,h_t}+t\,\Theta_{\det E,\det h_t}\otimes\Id_E>_P 0.\leqno(3.1)
$$
We start at $t=t_0$ and try to decrease $t$ as much as
we can, eventually going down to the positivity threshold $t_{\inf}=\tau_P(E)$.
In any case, we can find a $P$-positive metric on $E\otimes(\det E)^{t_0}$,
and since $\det(E\otimes(\det E)^{t_0})=(\det E)^{1+rt_0}$, we can derive
from it a metric $h_{t_0}$ on $E$ satisfying (3.1). In the sequel, we also set
$$
\omega_t=\Theta_{\det E,\det h_t}.\leqno(3.2)
$$
By our assumption (3.1), $\omega_t$ is a K\"ahler metric that lies in the
K\"ahler class $2\pi\,c_1(E)$. We wish to enforce suitable differential
equations on $(h_t)$ so that the family $(h_t)$ is uniquely determined,
running $t$ backwards as long as possible. One natural
condition is to require
$$
\Phi_P\big(\Theta_{E,h_t}+t\,\Theta_{\det E,\det h_t}\otimes\Id_E\big)
=\hbox{some positive volume form on $X$},
$$
in the hope of enforcing the $P$-positivity of the tensor
$\Theta_{E,h_t}+t\,\Theta_{\det E,\det h_t}\otimes\Id_E$.
Let $\Omega\in C^\infty(X,\Lambda^{n,n}T^*_X)$ be a fixed positive volume
form on $X$. For reasons that will become apparent later, we introduce
a new parameter $\beta\in\bR_+$ and the differential equation
$$
\Phi_P\big(\Theta_{E,h_t}+t\,\Theta_{\det E,\det h_t}\otimes\Id_E\big)
=f_t\,\bigg({\Omega\over\omega_t^n}\bigg)^{\!\beta}\,\Omega,
\quad f_t>0,~~f_t\in C^{\infty}(X,\bR).\leqno(3.3)
$$
In the case of Griffiths positivity, we
take $\Phi_P=\Phi_{G,s}$ with $p$ large, according to the discussion conducted
in 2.10~(e,f). As was pointed out in [Dem21], equality~(3.3) yields only
one scalar differential equation, whereas $h_t$ is represented by
$r^2$ unknown real coefficients.  Therefore we need to couple (3.3)
with an additional matrix equation of real rank $r^2-1$ to achieve exact
determinacy. It turns out that $r^2-1$  is precisely the real dimension
of trace free hermitian endomorphisms of~$E$. It is therefore natural
to consider trace free Hermite-Einstein equations of the form
$$
\omega_t^{n-1}\wedge\Theta^\circ_{E,h_t}=g_t\,\omega_t^n,\quad
g_t\in C^\infty(X,\Herm^\circ_{h_t}(E,E)),\leqno(3.3^\circ)
$$
expressed in terms of the direct sum decomposition
$$
\plainleqalignno{
\Herm^\circ_h(E,E)&=\big\{u\in\Herm_h(E,E)\,;\;\tr(u)=0\big\},&(3.4)\cr
\Herm_h(E,E)&=\Herm^\circ_h(E,E)\oplus\bR\,\Id_E,\quad
u=u^\circ+{1\over r}\tr(u)\otimes\Id_E,\quad \tr(u^\circ)=0.&(3.4')\cr}
$$
In the above notation, $\Theta_{E,h}^\circ$ is the curvature tensor of
$E\otimes (\det E)^{-1/r}$, namely
$$
\Theta_{E,h}^\circ=\Theta_{E,h}-
{1\over r}\Theta_{\det E,\det h}\otimes\Id_E\in
C^\infty(X,\Lambda_\bR^{1,1}T^*_X\otimes\Herm_h^\circ(E,E)).\leqno(3.5)
$$
By the fundamental work of [Don85] and [UhY86], we know that $(3.3^\circ)$
can be solved with $g_t=0$ if $E$ is $c_1(E)$-polystable, and with a
suitable choice of the right hand side \hbox{$g_t=G(h_t)$} otherwise;
as shown by [UhY86], it suffices to take for $G$ an appropriate
matrix functional, for instance $G(h)=-\varepsilon\log h$ in
suitable coordinates, with $\varepsilon>0$ arbitrary.
The ``friction term'' $g_t\,\omega_t^n=-\varepsilon\log h_t\,\omega_t^n$
helps in getting a priori bounds for the solutions, and in our case,
we will possibly need to take $\varepsilon$ large.
The~following simple observation is essential.
\medskip

\noindent
{\bf 3.6. Observation.} {\it As long as
$t\mapsto h_t$ is continuous with values in $C^2(X,\Hom(E,E))$ and we
start with an initial value $h_{t_0}$ such that
$\Theta_{E,h_t}+t\,\Theta_{\det E,\det h_t}\otimes\Id_E>_P0$ at time $t=t_0$,
complement {\rm (2.10~f)} shows that the positivity property $P$ is
preserved on the whole interval~$]t_{\inf},t_0]$ where the solution
exists. One would therefore need to show that the solution persists to
times $t<0$ to conclude that $0\in{}]t_{\inf},t_0]$ and
$(E,h_0)>_P0$.}  \medskip

\noindent
Now, the differential system $(3.2,\,3.3,\,3.3^\circ)$ can be considered with
arbitrary right hand sides $f_t$, $g_t$ of order at most $1$ in $h_t$,
i.e.\ of the form
$$
\plainleqalignno{
&f_t(z)=F(t,z,h_t(z),D_zh_t(z))>0,&(3.7)\cr
\noalign{\vskip6pt}
&g_t(z)=G(t,z,h_t(z),D_zh_t(z)),\quad
g_t\in C^\infty(X,\Herm^\circ_{h_t}(E,E)).&(3.7^\circ)
\cr}
$$
These right hand sides  do not affect the principal symbol of the system,
which is or order~$2$, as we will see very soon. At this stage, our
first concern is whether the above (fully non linear) differential system
is actually elliptic. It is a priori exactly determined in the sense
that there are exactly as many equations as unknowns,
namely $r^2$ scalar coefficients for $h_t(z)$.
\medskip

\noindent
{\bf 3.8. Theorem.} {\it Let $E\to X$ be a holomorphic vector bundle
such that $\det E$ is ample and \hbox{$t\in\bR$} such
that $E\otimes(\det E)^t>_P0$. Then there exist explicit distortion
functions $\beta_{P,t,h}$ in $C^0(X,\bR_+)$
such that for any metric $h_t$ on $E$ satisfying
$\Theta_{E,h_t}+t\,\Theta_{\det E,\det h_t}\otimes\Id_E>_P0$
and any \hbox{$\beta>\sup_X\beta_{P,h_t,t}$}, the system of differential
equations $(3.2,3.3,\,3.3^\circ)$ possesses an elliptic linearization in a
$C^2$ neighborhood of $h_t$, whatever is the choice of right hand sides
$f_t=F(t,z,h_t,D_zh_t)>0$, $g_t=G(t,z,h_t,D_zh_t)\in
\Herm^\circ_{h_t(z)}(E_z,E_z)$.}
\medskip

\noindent
{\it Proof.} The proof is similar to the one given in [Dem21], although
we have somewhat extended our perspective and allowed more flexible
equations. For simplicity of notation, we~put $h=h_t$ and, in general, we set
$$
\plaineqalign{
M&:=\Herm(E)=\hbox{hermitian forms $E\times E\to\bC$,~~
$M_+={}$positive ones in $M$},\cr
\noalign{\vskip4pt}
M_h&:=\Herm_h(E,E)=\hbox{hermitian endomorphisms $E\to E$ with respect
to $h\in M_+$},\cr
\noalign{\vskip4pt}
M_h^\circ&:=\Herm_h^\circ(E,E)
=\hbox{trace free hermitian endomorphisms $E\to E$.}\cr}
$$
The system of equations $(3.2,3.3,3.3^\circ)$ is associated with the nonlinear
differential operator
$$
Q:C^\infty(X,M_+)\to C^\infty(X,\bR\oplus M_h^\circ),\qquad
h\mapsto Q(h)
$$
defined by $Q=Q_\bR\oplus Q^\circ$ where
$$
\plaincases{
\omega_h:=\Theta_{\det E,\det h}>0,\phantom{\big.}\cr
\noalign{\vskip8pt}
Q_\bR(h):=\big(\omega_h^n/\Omega)^\beta\,\Omega^{-1}\,\Phi_P\big(\Theta_{E,h}
+t\,\omega_h\otimes\Id_E\big),\cr
\noalign{\vskip8pt}
Q^\circ(h):=(\omega_h^n)^{-1}\big(\omega_h^{n-1}\wedge\Theta_{E,h}^\circ\big).\cr}
\leqno(3.9)
$$
By definition, it is elliptic at $h$ 
if its linearization $dQ(h)$ is an elliptic linear operator,
the exact determinacy being reflected in the fact that $M$ and
\hbox{$\bR\oplus M_h^\circ$} have the same rank $r^2$ over the field~$\bR$
of real numbers. Our goal is to compute the principal symbol
$$
\sigma_2(dQ(h))\in C^\infty(X,S^2T^\bR_X\otimes\Hom(M,\bR\oplus M_h^\circ))
$$
of the linearized operator $dQ(h)$, and to check that
$\sigma_2(dQ(h))(\xi)\in\Hom(M,\bR\oplus M_h^\circ)$ is invertible
for every non zero vector cotangent vector $\xi\in T^*_X$.
For the calculation in coordinates, we fix locally on~$X$ a holomorphic frame
$(\varepsilon^0_\lambda)_{1\leq \lambda\leq r}$ of $E$, and denote by $H_0$ the
trivial hermitian metric for which $(\varepsilon^0_\lambda)$ is orthonormal.
Any hermitian metric $h$ is then represented by a hermitian matrix,
denoted again $h=(h_{\lambda\mu})$, such that the corresponding inner
product is $\langle h\bu,\bu\rangle_{H_0}$. It is well known that the
Chern curvature tensor $\Theta_{E,h}$ is given locally by the matrix of
$(1,1)$-forms
$$
\Theta_{E,h}=i\ol\partial(h^{-1}\partial h).
$$
Next, we pick an infinitesimal variation
$\delta h$ of $h$ in $C^\infty(X,M)$. It is convenient to write
it under the form $\delta h=\langle u\,\bu,\bu\rangle_h=
\langle hu\,\bu,\bu\rangle_{H_0}$ with $u\in M_h=\Herm_h(E,E)$.
In terms of matrices in a fixed $H_0$-orthonormal frame,
we thus have $\delta h=hu$, i.e., $u=(u_{\lambda\mu})=h^{-1}\delta h$ is
some sort of ``logarithmic variation of~$h$''. In~this setting, we first
evaluate $(d\Theta_{E,h})(\delta h)$. We have $h+\delta h=h(I+u)$ and
$(h+\delta h)^{-1}=(\Id-u)h^{-1}$ modulo~$O(u^2)$. As an abuse
of notation, we will write $d\Theta_{E,h}(u)$ what we should write
$d\Theta_{E,h}(\delta h)$, and likewise for the other differentials;
of course this makes no difference in terms of matrices if at the end
of the calculation we express the result in an $h$-orthonormal (rather
than $H_0$-orthonrmal) frame of $E$, in which $h$ is represented by
the unit matrix. We obtain
$$
d\Theta_{E,h}(u)
=i\ol\partial\big((I-u)h^{-1}\partial(h(I+u)\big)-
i\ol\partial(h^{-1}\partial h)\quad{\rm mod}~O(|u|^2+|du|^2),
$$
that is,
$$
d\Theta_{E,h}(u)
=\ol\partial(h^{-1}\partial(hu))
-i\ol\partial(uh^{-1}\partial h)
=-\partial\ol\partial u+\ol\partial((h^{-1}\partial h)u)
-i\ol\partial(u(h^{-1}\partial h)).
\leqno(3.10)
$$
As a consequence, the order~$2$ part $(\bu)^{[2]}$ of the linearized
operator $d(\Theta_{E,h})$, in other words its principal symbol is
simply given by
$$
(d\Theta_{E,h})^{[2]}(u)=-i\partial\ol\partial u.\leqno(3.10')
$$
Since $d\omega_h(u)=\tr(d\Theta_{E,h}(u))$, we find
$$
\plainleqalignno{
&(d\omega_h)^{[2]}(u)=-i\tr\partial\ol\partial u=-i\partial\ol\partial\tr(u),
\qquad (d\Theta_{E,h}^\circ)^{[2]}(u)=-i\partial\ol\partial u^\circ,&(3.11)\cr
\noalign{\vskip4pt}
&(d\log\omega_h^n)^{[2]}(u)=
n\,(\omega_h^n)^{-1}\omega_h^{n-1}\wedge (d\omega_h)^{[2]}(u).&(3.11')\cr}
$$
In order to compute $dQ$, we need the differential of the functional $\Phi_P$.
In the case $P=N,N^*$, we have to consider the $(nr\times nr)$-matrix
$\theta_t(h)$ of the hermitian form on $T_X\otimes E$ defined by
$$
\theta_t(h)\simeq\Theta_{E,h}+t\,\omega_h\otimes\Id_E>_P0
$$
(or its transpose), and the logarithmic differential of
$\det(\theta_t(h))^{1/r}$ is ${1\over r}\,\tr\big(\theta_t(h)^{-1}d\theta_t(h)
\big)$ where $\theta_t(h)^{-1}=(\det\theta_t(h))^{-1}{\;}^T(\theta_t(h)^{\rm cof})$
and $\theta_t(h)^{\rm cof}$ is the $\Hom(T_X\otimes E,T_X\otimes E)$-cofactor
matrix of $\theta_t(h)$, ${}^T(\bu)$ the corresponding transposition
operator. We pursue our calculations with respect to $\omega_h$-orthonormal
coordinates $(z_j)_{1\le j\le n}$ on~$X$ at a given point $z^0\in X$,
and also use later on an $h$-orthonormal frame
$(e_\lambda)_{1\leq\lambda \leq r}$ of~$E_{z_0}$. We then get respectively
$$
\plainleqalignno{\kern43pt
&d\theta_t(h)^{[2]}(u)=
-\bigg({\partial^2 u_{\lambda\mu}\over\partial z_j\partial\ol z_k}
+t\,\delta_{\lambda\mu}
\sum_\nu{\partial^2 u_{\nu\nu}\over\partial z_j\partial\ol z_k}
\bigg)_{(j,\lambda),(k,\mu)},&(3.12)\cr
&(d\log\omega_h^n)^{[2]}(u)=-\sum_{j,\nu}
{\partial^2 u_{\nu\nu}\over\partial z_j\partial\ol z_j},&(3.13)\cr
&d\log\Phi_N(\theta_t)_h^{[2]}(u)={-1\over r\,\det\theta_t(h)}
\sum_{j,k,\lambda,\mu}\theta_t(h)^{\strut\rm cof}_{jk\lambda\mu}
\bigg({\partial^2 u_{\lambda\mu}\over\partial z_j\partial\ol z_k}
+t\,\delta_{\lambda\mu}
\sum_\nu{\partial^2 u_{\nu\nu}\over\partial z_j\partial\ol z_k}
\bigg),&(3.14_N)\cr
&d\log\Phi_{N^*}(\theta_t)_h^{[2]}(u)={-1\over r\,\det{}^T\theta_t(h)}
\sum_{j,k,\lambda,\mu}\!\big({\,}^T\theta_t(h)\big)^{\rm cof}_{jk\mu\lambda}
\bigg({\partial^2 u_{\lambda\mu}\over\partial z_j\partial\ol z_k}
+t\,\delta_{\lambda\mu}
\sum_\nu{\partial^2 u_{\nu\nu}\over\partial z_j\partial\ol z_k}
\bigg)&(3.14_{N^*})
\cr}
$$
where $\big({}^T\theta_t(h)\big)^{\rm cof}$ is the
$\Hom(T_X\otimes E^*,T_X\otimes E^*)$-cofactor matrix of ${}^T\theta_t(h)$.
The~calculation for the functional $\Phi_{G,s}$ requires a differentiation of
(2.17) and is more involved. If we notice that the differentiation of
$\langle\bu,\bu\rangle_h$ in $h$ does not contribute to the order 2 terms,
we find
$$
\plainleqalignno{
&\kern44pt d\log\Phi_{G,s}(\theta_t)_h^{[2]}(u)=&(3.14_G)\cr
&\Bigg(\int_{\textstyle{v\in E\atop |v|_h=1}}
{d\sigma(v)\over
\big((\langle\theta_t(h)\cdot v,v\rangle_h)^n\big)^{s\phantom{!\!}}}
\Bigg)^{\!-1}
\int_{\textstyle{v\in E\atop|v|_h=1}}
{n\,(\langle\theta_t(h)\cdot v,v\rangle_h)^{n-1}\wedge
\langle d\theta_t(h)^{[2]}(u)\cdot v,v\rangle_h\,d\sigma(v)
\over\big((\langle\theta_t(h)\cdot v,v\rangle_h)^n\big)^{s+1\phantom{!\!}}}.\cr
}
$$
Notice that the $s$-th power of an $(n,n)$-form in the first integral and
the quotient of an $(n,n)$-form by the $(s+1)$-st power of an $(n,n)$-form
in the second integral actually combine into a dimensionless value.
In normal coordinates,
$\langle d\theta_t(h)^{[2]}(u)\cdot v,v\rangle_h$ is the $(1,1)$-form
$$
\langle d\theta_t(h)^{[2]}(u)\cdot v,v\rangle_h=
-\sum_{j,k,\lambda,\mu}\bigg(
{\partial^2 u_{\lambda\mu}\over\partial z_j\partial\ol z_k}\,v_\lambda
\overline v_\mu+t\,{\partial^2 u_{\lambda\lambda}\over\partial z_j\partial\ol z_k}
\,|v_\mu|^2\bigg)\,dz_j\wedge d\overline z_k.
\leqno(3.15)
$$
Let us begin with the case of Nakano positivity $P=N$. By the above identities,
the logarithmic differential of the first scalar component $Q_\bR(h)$
of $Q(h)$ has order 2 terms
$$
\plainleqalignno{
\kern2pt&Q_{\bR}(h)^{-1}\,(dQ_{\bR,h})^{[2]}(u)
=d\log\Phi_N(\theta_t)_h^{[2]}(u)+\beta\,(d\log\omega_h^n)^{[2]}(u)\cr
&\kern34pt{}={-1\over r\,\det\theta_t(h)}
\sum_{j,k,\lambda,\mu}\theta_t(h)^{\strut\rm cof}_{jk\lambda\mu}
\bigg({\partial^2 u_{\lambda\mu}\over\partial z_j\partial\ol z_k}
+t\,\delta_{\lambda\mu}
\sum_\nu{\partial^2 u_{\nu\nu}\over\partial z_j\partial\ol z_k}\bigg)
-\beta\sum_{j,\nu}{\partial^2 u_{\nu\nu}\over\partial z_j\partial\ol z_j},\cr
&\kern34pt{}={-1\over r\,\det\theta_t(h)}
\sum_{j,k,\lambda,\mu}\theta_t(h)^{\strut\rm cof}_{jk\lambda\mu}
\bigg({\partial^2 u^\circ_{\lambda\mu}\over\partial z_j\partial\ol z_k}
+\Big(t+{1\over r}\Big)\,\delta_{\lambda\mu}
\sum_\nu{\partial^2 u_{\nu\nu}\over\partial z_j\partial\ol z_k}\bigg)
-\beta\sum_{j,\nu}{\partial^2 u_{\nu\nu}\over\partial z_j\partial\ol z_j},
&(3.16_N)\cr}
$$
and we get a similar expression for $\Phi_{N^*}$ by $(3.14_{N^*})$.
Finally, we compute the order 2 terms in the differential of the
second component
$$
h\mapsto Q^\circ(h)=(\omega_h^n)^{-1}\,(\omega_h^{n-1}\wedge\Theta^\circ_{E,h})
\in M^\circ_h.
$$
The above calculations imply
$$
\plaineqalign{
dQ^\circ(h)^{[2]}(u)
={}&-n(\omega_h^n)^{-2}\,\big(\omega_h^{n-1}\wedge(d\omega_h)^{[2]}(u)\big)\cdot
\big(\omega_h^{n-1}\wedge \Theta^\circ_{E,h}\big)\cr
&+(n-1)\,(\omega_h^n)^{-1}\,\big(\omega_h^{n-2}\wedge (d\omega_h)^{[2]}(u)
\wedge\Theta^\circ_{E,h}\big)\cr
&+(\omega_h^n)^{-1}\,\big(\omega_h^{n-1}\wedge(d\Theta^\circ_{E,h})^{[2]}(u)
\big).\cr}
$$
If we denote $\Theta^\circ_{E,h}=\sum_{j,k,\lambda,\mu}
c^\circ_{jk\lambda\mu}\,dz_j\wedge d\overline z_k\otimes e^*_\lambda\otimes
e_\mu$ at $z^0$, this yields
$$
\plainleqalignno{
dQ^\circ(h)^{[2]}(u)
={}&+{1\over n}\sum_{j,\nu}{\partial^2 u_{\nu\nu}\over\partial z_j\partial\ol z_j}
\cdot\sum_{k,\lambda,\mu}c^\circ_{kk\lambda\mu}\,e_\lambda^*\otimes e_\mu\cr
&-{1\over n}\sum_{j,k,\lambda,\mu,\nu}\bigg(
{\partial^2 u_{\nu\nu}\over\partial z_j\partial\ol z_j}\,
c^\circ_{kk\lambda\mu}\,e_\lambda^*\otimes e_\mu-
{\partial^2 u_{\nu\nu}\over\partial z_k\partial\ol z_j}\,
c^\circ_{jk\lambda\mu}\,e_\lambda^*\otimes e_\mu\bigg)\cr
&-{1\over n}\sum_{j,\lambda,\mu}
{\partial^2 u^\circ_{\lambda\mu}\over\partial z_j\partial\ol z_j}
\,e_\lambda^*\otimes e_\mu\cr
={}&{1\over n}\sum_{j,k,\lambda,\mu,\nu}
{\partial^2 u_{\nu\nu}\over\partial z_k\partial\ol z_j}\,
c^\circ_{jk\lambda\mu}\,e_\lambda^*\otimes e_\mu-{1\over n}\sum_{j,\lambda,\mu}
{\partial^2 u^\circ_{\lambda\mu}\over\partial z_j\partial\ol z_j}
\,e_\lambda^*\otimes e_\mu.&(3.16^\circ)\cr}
$$
The principal symbol $\sigma_2(dQ(h))$ at $h$, taken on a cotangent vector
$\xi\in T^*_X$, is thus given by the two components $\sigma_2(dQ_\bR(h))$
and $\sigma_2(dQ^\circ(h))$ such that
$$
\plainleqalignno{
{\sigma_2(dQ_\bR(h))(\xi)\cdot u\over dQ_\bR(h)}
&={-1\over r\,\det\theta_t(h)}
\sum_{j,k,\lambda,\mu}\!\theta_t(h)^{\strut\rm cof}_{jk\lambda\mu}\,
\xi_j\overline\xi_k\bigg(u^\circ_{\lambda\mu}+\Big(t+{1\over r}\Big)\,
\delta_{\lambda\mu}\tr(u)\bigg)-\beta\,|\xi|^2\tr(u),\cr
&&(3.17,3.17^\circ)\cr
\noalign{\vskip-4pt}
\sigma_2(dQ^\circ(h))(\xi)\cdot u
&=-{1\over n}\,{\omega_h^n\over\Omega}\,\sum_{\lambda,\mu}
\Bigg(\sum_{j\neq k}\big(c^\circ_{kk\lambda\mu}\,|\xi_j|^2-c^\circ_{jk\lambda\mu}
\xi_j\overline\xi_k\big)\tr(u)
+ |\xi|^2\,u^\circ_{\lambda\mu}\Bigg)\,e_\lambda^*\otimes e_\mu.\cr}
$$
By definition $dQ(h)$ is elliptic if and only if 
$\sigma_2(dQ(h))(\xi)\in\Hom(M_h,\bR\oplus M^\circ_h)$
is injective for all cotangent vectors $\xi\neq 0$. Now, since
$t+{1\over r}>0$ and since the cofactor matrix is hermitian
positive by the Nakano positivity assumption, we see that the
vanishing of $\sigma_2(dQ_\bR(h))(\xi)\cdot u$ implies
by a simple Cauchy-Schwarz argument that
$$
|\tr(u)|\leq {1\over\beta r}\,{|\theta_t(h)^{\rm cof}|\over \det\theta_t(h)}\,
|u^\circ|,\leqno(3.18)
$$
where the norms of tensors are taken with respect to $(E,h)$
and $(T_X,\omega_h)$. By plugging this inequality into
$\sigma_2(dQ^\circ(h))$, taking the inner product with
$u^\circ=\sum_{\lambda,\mu}u^\circ_{\lambda\mu}\,e^*_\lambda\otimes e_\mu$
and using again Cauchy-Schwarz,
we see that $\sigma_2(dQ_\bR(h))(\xi)\cdot u=0$ entails
$$
\langle-\sigma_2(dQ^\circ(h))(\xi)\cdot u,u^\circ\rangle\geq
{1\over n}\,{\omega_h^n\over\Omega}\,\Big(|\xi|^2\,|u^\circ|^2-
(\sqrt{n-1}+1)\,|\Theta^\circ_{E,h}|\,|\xi|^2\,|u^\circ|\,|\tr(u)|\Big),
$$
hence
$$
\big|\sigma_2(dQ^\circ(h))(\xi)\cdot u\big|\geq
{1\over n}\,{\omega_h^n\over\Omega}\,|\xi|^2\,|u^\circ|\,
\bigg(1-{\sqrt{n-1}+1\over\beta r}\,{|\Theta^\circ_{E,h}|\,|\theta_t(h)^{\rm cof}|
\over \det\theta_t(h)}\bigg).
\leqno(3.18^\circ)
$$
Let us introduce the ``distortion function'' $\beta_{N,t,h}\in C^0(X,\bR_+)$
$$
\beta_{N,t,h}={\sqrt{n-1}+1\over r}\,{|\Theta^\circ_{E,h}|\,
|\theta_t(h)^{\rm cof}|\over \det\theta_t(h)}\leqno(3.19_N)
$$
computed at each point $z\in X$ and, in a similar manner for dual
Nakano positivity,
$$
\beta_{N^*,t,h}={\sqrt{n-1}+1\over r}\,
{|\Theta^\circ_{E,h}|\,|({\kern1pt}^T\theta_t(h))^{\rm cof}|
\over \det({}^T\theta_t(h))}.\leqno(3.19_{N^*})
$$
Then, for $\beta>\sup_X\beta_{N,t,h}$ (resp.\
$\beta>\sup_X\beta_{N^*,t,h}$), inequalities
(3.18) and $(3.18^\circ)$ imply the ellipticity of our differential system.
In the case of $\Phi_{G,s}$, the identities $(3.14_G)$ and (3.15) yield
$$
\plaineqalign{
&Q_{\bR}(h)^{-1}\,\sigma_2(dQ_{\bR,h})(u)\cdot\xi
=\sigma_2(d\log\Phi_{G,s}(\theta_t)_h)(u)\cdot\xi+
\beta\,\sigma_2(d\log\omega_h^n)(u)\cdot\xi\cr
&\kern20pt{}=-\beta\,|\xi|^2\tr(u)-\Bigg(\int_{\textstyle{v\in E\atop |v|_h=1}}
{d\sigma(v)\over
\big((\langle\theta_t(h)\cdot v,v\rangle_h)^n\big)^{s\phantom{!\!}}}
\Bigg)^{\!-1}\cr
&\kern40pt{}\times\int_{\textstyle{v\in E\atop|v|_h=1}}
{n\,(\langle\theta_t(h)\cdot v,v\rangle_h)^{n-1}\wedge\big(
\langle u^\circ(v),v\rangle+(t+{1\over r})\tr(u)\,|v|^2\big)\,
i\xi\wedge\overline\xi\,d\sigma(v)
\over\big((\langle\theta_t(h)\cdot v,v\rangle_h)^n\big)^{s+1\phantom{!\!}}}.\cr}
$$
By easy estimates left to the reader, this leads in the case of the
Griffiths functional to the distortion function
$$
\plainleqalignno{\beta_{G,s,t,h}&=(\sqrt{n-1}+1)\,|\Theta^\circ_{E,h}|\cr
&\kern20pt{}\times\Bigg(\int_{\textstyle{v\in E\atop |v|_h=1}}{d\sigma(v)\over
\big((\langle\theta_t(h)\cdot v,v\rangle_h)^n\big)^{s\phantom{!\!}}}
\Bigg)^{\!-1}\int_{\textstyle{v\in E\atop|v|_h=1}}
{n\,(\langle\theta_t(h)\cdot v,v\rangle_h)^{n-1}\wedge \omega_h\,d\sigma(v)
\over\big((\langle\theta_t(h)\cdot v,v\rangle_h)^n\big)^{s+1\phantom{!\!}}},
&(3.19_G)\cr}
$$
where $\beta_{G,s,t,h}(z)$ is obtained by computing the integrals
fiberwise on~$E_z$. The proof of Theorem~3.8 is complete.\hfill\qed
\medskip

\noindent
{\bf 3.20. Remark.} It the curvature tensor $\Theta_{E,h}(z)$
happens to be just rescaled by a positive multiplication factor
at some point $z\in X$, the value of the above distortion
function $\beta_{P,t,h}(z)$ can be seen to remain invariant. In some sense,
$\beta_{P,t,h}(z)$ measures the ratio of ``eigenvalues'' along directions of
maximum and minimum $P$-positivity for $\Theta_{E,h}+t\,\Theta_{\det E,\det h}
\otimes\Id_E$, at each point $z\in X$. Hopefully, there might be a way
of relating these distortion functions to simpler geometric invariants,
such as the slopes in the Harder-Narasimhan filtration of $E$ with 
respect to $c_1(E)$.
\medskip

Our next concern is to ensure that the existence and uniqueness of
solutions hold, at least on suitable subsets of
$\bR\times C^\infty(X,M_+)$, consisting of pairs $(t,h)$ such that
$\theta_{E,h}+t\,\Theta_{\det E,\det h}\otimes\Id_E>_P0$. We fix such a pair
$(t_0,h_0)$ and use $H_0=h_{t_0}$ and $\omega_{t_0}=\Theta_{\det E,\det h_{t_0}}$
as the reference metrics on $E$ and $T_X$ respectively.
For $K\geq K_0\gg 1$, we consider the subset
$S_K\subset{}]-1/r,t_0]\times C^\infty(X,M_+)$ of pairs $(t,h)$
such that
$$
|h|_{C^2}\leq K,\quad |h^{-1}|_{C^2}\leq K,\quad
\theta_{E,h}+t\,\Theta_{\det E,\det h}\otimes\Id_E
\geq_P K^{-1}\,\omega_0\otimes\Id_E\leqno(3.21)
$$
with respect to $(H_0,\omega_{t_0})$. In~the case of a rank one metric
$h=e^{-\varphi}$, it is well-known that the K\"ahler-Einstein equation
\hbox{$(\omega_0+i\partial\ol\partial\varphi_t)^n=e^{tf+\lambda\varphi_t}
\omega_0^n$} yields easily the openness and closedness of solutions
when $\lambda$ is positive, as a consequence of the fact that the
linearized operator $\psi\mapsto \Delta_{\omega_{\varphi_t}}\psi-\lambda \psi$
is~always invertible. Here we can still play the game of
adjusting the right hand sides $f_t,g_t$ in $(3,3,3.3^\circ)$ to
achieve the invertibility of the related elliptic operator $\widehat Q$,
at least for $(t,h)\in S_K$. Before doing so, we introduce
some notation. If $h\in\Herm(E)$ is a hermitian form, we have an isomorphism
$$
\Herm(E)\to\Herm_{h_{t_0}}(E,E),\quad
h\mapsto\tilde h~~\hbox{such that}~~
h(v,w)=\langle v,w\rangle_h=\langle\tilde h(v),w\rangle_{h_{t_0}},\leqno(3.22)
$$
and for $h\in\Herm_+(E)$, we let $\log\tilde h\in\Herm_{h_{t_0}}(E,E)$ be
its logarithm as a hermitian endomorphism. Finally, we define
$\tilde h^{(1)}=(\det\tilde h)^{-1/r}\tilde h$, so that
$\det(\tilde h^{(1)})=1$ and
$$
\log\tilde h^{(1)}=(\log\tilde h)^\circ~\in~\Herm^\circ_{h_{t_0}}(E,E)
\leqno(3.23)
$$
is the trace free part of $\log\tilde h$.
One way to generalize the K\"ahler-Einstein condition to the case of
arbitrary ranks $r\ge 1$ is to consider pairs $(t,h)$ satisfying
a differential equation of the form
(3.3) with a factor $f_t(z)=(\det h_{t_0}(z)/\det h(z))^\lambda$, namely
$$
\Phi_P\big(\Theta_{E,h}+t\,\omega_h\otimes\Id_E\big)
=\bigg({\det h_{t_0}\over\det h}\bigg)^{\!\lambda}
\,\bigg({\Omega\over\omega_h^n}\bigg)^{\!\beta}\,\Omega,\quad
\hbox{where}~~\omega_h=\Theta_{\det E,\det h},~\lambda,\beta>0,
\leqno(3.24)
$$
and the volume form $\Omega>0$ is chosen so that equation (3.24) is
satisfied by $(t_0,h_{t_0})$. The choice $\lambda>0$ has the advantage that
the right hand side gets automatically rescaled when multiplying $h$
by a constant (while the left hand side remains untouched), thus avoiding
a trivial non invertibility issue. When $r=1$, one easily sees that
equation (3.24) actually reduces to the usual K\"ahler-Einstein equation.
By Uhlenbeck-Yau [UhY86], if one chooses for the right hand side
of $(3.3^\circ)$ a ``friction term'' $g_t$ of the type
$g_t(z)=-\varepsilon\,a(t,z)\,\log\tilde h^{(1)}(z)$, $a(t,z)>0$, then
the Hermite-Einstein equation always has a solution, although it
usually blows up as $\varepsilon\to 0$ when $E$ is unstable.
This leads to couple (3.24) with a trace free Hermite-Einstein
equation of the form
$$
\omega_h^{-n}\,\big(\omega_h^{n-1}\wedge\Theta^\circ_{E,h}\big)=
-\varepsilon\,A(\det h)\,\log\tilde h^{(1)},
\leqno(3.24^\circ)
$$
where $A\in C^\infty(]0,+\infty[\,,\bR_+)$ is a positive function;
one could also use more generally a factor $A(\det u(z),z)$ where
$A\in C^\infty(\,]0,+\infty[{}\times X,\bR_+)$, such as 
$A(y,z)=(\det h_{t_0}(z)/y)^\mu$, $\mu\in\bR$. The precise form of $A$
is irrelevant here, provided that $A>0$; one could just take $A\equiv 1$.
The~right hand sides used in $(3.24,3.24^\circ)$ do not depend on
higher derivatives of $h$, thus Theorem~3.8 ensures the ellipticity
of the differential system as soon $\beta>\sup_X\beta_{P,t,h}$
(see~$(3.19_P)$).
\medskip

\noindent
{\bf 3.25. Theorem.} {\it Consider the differential operator
$\widehat Q:C^\infty(X,M_+)\to C^\infty(X,\bR\oplus M_h^\circ)$
defined by
$$
\plainleqalignno{
\widehat Q_\bR(h)
&=\bigg({\det h\over\det h_{t_0}}\bigg)^{\!\lambda}\,Q_\bR(h)
=\bigg({\det h\over\det h_{t_0}}\bigg)^{\!\lambda}\,
\bigg({\omega_h^ n\over\Omega}\bigg)^{\!\beta}\,\Omega^{-1}\,
\Phi_P\big(\Theta_{E,h}+t\,\omega_h\otimes\Id_E\big),&({\rm YM})\cr
\widehat Q^\circ(h)&=Q^\circ(h)+\varepsilon\,A(h)\,\log\tilde h^{(1)}
=\omega_h^{-n}\big(\omega_h^{n-1}\wedge\Theta^\circ_{E,h}\big)+\varepsilon\,
A(\det h)\,\log\tilde h^{(1)}.
&({\rm YM}^\circ)\cr}
$$
where $Q=Q_\bR\oplus Q^\circ$ is the operator introduced in the
proof of Theorem~$3.8$, and $A$ is any smooth positive function.
There exist bounds $\beta_0(K):=\sup_{(t,h)\in S_K}\sup_X\beta_{P,t,h}$,
$\varepsilon_0(A,K,\beta)$ and $\lambda_0(A, K,\beta)$ such
that for any choice of cons\-tants
$\beta>\beta_0(K)$, $\varepsilon>\varepsilon_0(A,K,\beta)$ and
$\lambda>\lambda_0(A,K,\beta)$, the elliptic operator defined by
$({\rm YM},~{\rm YM}^\circ)$ possesses an invertible elliptic linearization
$\smash{d\widehat Q(h)}$ for all $(t,h)\in S_K$. As a consequence,
there exists an open interval $[t_0-\delta_0,t_0]$, $\delta_0>0$,
such that the solution $h_t$ of the system
$\smash{(\widehat Q_\bR(h),\widehat Q^\circ(h))=(1,0)}$ exists
and is unique for $t\in[t_0-\delta_0,t_0]$. This solution $h_t$ depends
differentiably on~$t$.}\medskip

\noindent
{\it Proof.} Here, we have to keep an eye on the linearized operator
$d\widehat Q$ itself, and not just its principal symbol. We let again
$u=h^{-1}\delta h\in \Herm_h(E,E)$ and use the formulas established
for $dQ(h)$ in the proof of Theorem~3.8. The logarithmic derivative
of $\smash{\widehat Q_{\bR}(h)}$ is
$$
\widehat Q_{\bR}(h)^{-1}\,d\widehat Q_{\bR}(h)(u)
=Q_{\bR}(h)^{-1}\,dQ_{\bR}(h)(u)+\lambda\,\tr(u).\leqno(3.26)
$$
For $\widehat Q^\circ$, we need the fact that, when viewed as a
hermitian endomorphism, $h^\circ=h\cdot(\det h)^{-1/r}$ possesses
a~logarithmic variation
$$
(\tilde h^{(1)})^{-1}\delta \tilde h^{(1)}= u^\circ
= u -{1\over r}\tr(u)\cdot\Id_E.
$$
By the classical formula expressing the differential of
the logarithm of a matrix, we have
$$
(d\log g)(\delta g)=
\int_0^1\big((1-s)\Id+sg\big)^{-1}\delta g\,\big((1-s)\Id+sg\big)^{-1}\,ds
$$
($g$ and $\delta g$ need not commute here!), which implies
$$
d\log\tilde h^{(1)}(\delta h)=
\int_0^1\big((1-s)\Id+s\,\tilde h^{(1)}\big)^{-1}
\,\tilde h^{(1)} u^\circ\,\big((1-s)\Id+s\,\tilde h^{(1)}\big)^{-1}\,ds.
$$
If $(\alpha_\lambda)_{1\leq\lambda\leq r}$ are the eigenvalues of $h^{(1)}$
with respect to $h_{t_0}$ and we use an orthonormal basis of
eigenvectors, we obtain in coordinates
$$
d\log\tilde h^{(1)}:
\delta h\longmapsto L_h(u^\circ)=\big(\gamma_{\lambda\mu}u^\circ_{\lambda\mu}
\big)_{1\leq\lambda,\mu\leq r},\quad
\gamma_{\lambda\mu}={\alpha_\mu\over \alpha_\mu-\alpha_\lambda}
\log{\alpha_\mu\over\alpha_\lambda},
$$
where $u=(u_{\lambda\mu})_{1\leq\lambda,\mu\leq r}=h^{-1}\delta h$ and the
coefficient $\gamma_{\lambda\mu}>0$ is to be interpreted as~$1$ if
$\alpha_\lambda=\alpha_\mu$. In the end, we obtain
$$
d\widehat Q^\circ(h)(u)=dQ^\circ(h)(u)+
\varepsilon\,\big(A(h)\,L_h(u^\circ)+A'(\det h)\,
\det h\, \tr(u)\,\log\tilde h^{(1)}\big).
\leqno(3.26^\circ)
$$
In order to check the invertibility, we compare the operators
$$
d\log Q_\bR(h)\oplus dQ^\circ(h)\quad\hbox{and}\quad
d\log\widehat Q_\bR(h)\oplus d\widehat Q^\circ(h).
$$
The principal symbol calculations $(3.17,3.17^\circ)$ 
show that for $\beta>\beta_0(K)=\sup_{(t,h)\in S_k}\sup_X\beta_{P,h,t}$,
the linearized operator $d\log Q_\bR(h)\oplus dQ^\circ(h)$ is
elliptic and essentially positive for all $(t,h)\in S_K$.
We consider the natural $L^2$ metric on $L^2(X,M)\simeq
L^2(X,\bR\oplus M_h^\circ)$ defined by
$$
\Vert u\Vert^2=\Vert\tr(u)\Vert^2+\Vert u^\circ\Vert^2,
$$
using the hermitian metric $h^*\otimes h$ on $\Herm_h(E,E)$ and the
volume element $\omega_h^n/n!$ on $X$. By the ellipticity of $dQ(h)$ and
an elementary case of G{\aa}rding's inequality, there exist
constants $C,C'>0$ such that
$$
\langle\!\langle d\log Q_\bR(h)(u)\oplus dQ^\circ(h)(u),
\tr(u)\oplus u^\circ \rangle\!\rangle_{\bR\oplus M_h^\circ}\geq
C\,\Vert\nabla_hu\Vert^2-C'\Vert u\Vert^2.\leqno(3.27)
$$
Moreover, the estimate is valid with constants $C=C(K,\beta)$,
$C'=C'(K,\beta)$, uniformly for all $(t,h)\in S_K$. In such a $C^2$ bounded set,
we also have bounds
$$
\plaincases{
\langle A(h)\,L_h(u^\circ),u^\circ\rangle =
A(h)\sum_{1\leq\lambda,\mu}\gamma_{\lambda\mu}|u^\circ_{\lambda\mu}|^2
\geq C''\,|u^\circ|^2,\cr
\noalign{\vskip4pt}
|A'(\det h)\,\det h\,\log\tilde h^{(1)}|\leq C'''\cr}
\leqno(3.28)
$$
with $C''=C''(A,K),C'''=C'''(A,K)>0$.
Estimates $(3.26,3.26^\circ)$, $(3.27,3.28)$
and the Cauchy-Schwarz inequality imply
$$
\plainleqalignno{
&\langle\!\langle d\log\widehat Q_\bR(h)(u)\oplus d\widehat Q^\circ(h)(u),
\tr(u)\oplus u^\circ \rangle\!\rangle_{\bR\oplus M_h^\circ}\cr
\noalign{\vskip3pt}
&\kern40pt{}\geq
C\,\Vert\nabla_hu\Vert^2-C'\Vert u\Vert^2
+\lambda\Vert\tr(u)\Vert^2+\varepsilon\,\big(
C''\,\Vert u^\circ\Vert^2-C'''\,
\Vert\tr(u)\Vert\,\Vert u ^\circ\Vert\big)\cr
&\kern40pt{}\geq
C\,\Vert\nabla_hu\Vert^2+\bigg(\lambda-C'-{2(C''')^2\over 2C''}\bigg)\,
\Vert\tr(u)\Vert^2+\bigg({1\over 2}\,\varepsilon\,C''-C'\bigg)
\Vert u ^\circ\Vert^2&(3.29)\cr}
$$
by the inequality $\Vert\tr(u)\Vert\,\Vert u ^\circ\Vert\leq
{C'''\over 2C''}\Vert\tr(u)\Vert^2+{C''\over 2C'''}\Vert u ^\circ\Vert^2$.
If we take
$$
\varepsilon>\varepsilon_0(A,K,\beta)={2C'\over C''},\quad
\lambda>\lambda_0(A,K,\beta)=C'+{2(C''')^2\over 2C''},\leqno(3.30)
$$
we conclude from (3.29) that
$d\log\widehat Q_\bR(h)\oplus d\widehat Q^\circ(h)$ is an invertible
elliptic operator $W^{s+2}\to W^s$ for all Sobolev spaces $W^s$,
$s\geq 0$. The proof of Theorem~3.25 is achieved by applying standard
results in the theory of elliptic operators and the implicit function
theorem.\hfill\qed
\medskip

\noindent
{\bf 3.31. Remark.} Theorem 3.25 is somehow purely local. The main point
would be to obtain more uniform estimates with respect to the
metric $h$, especially in terms of the distortion fonctions, so that one
could keep control on the solution throughout the expected maximal interval
of time. This obviously requires a finer analysis than the one we
conducted here. If $\beta_{P,t,h}$ could be better understood, 
explicit expressions of the constants $C,C',C'',C'''$
and thus of $\varepsilon_0(A,K,\beta)$ and $\lambda_0(A,K,\beta)$ would 
perhaps become accessible by looking more in depth at the Bochner formula.
\bigskip\bigskip

\centerline{\bigbf References}
\vskip8pt

\Bibitem[Ber09]&Berndtsson B.:& Curvature of vector bundles associated to
holomorphic fibrations,& Annals of Math.\ {\bf 169} (2009), 531–-560&

\Bibitem [DemB]&Demailly, J.-P.:& Complex analytic and algebraic geometry,&
Online book at\hfil\break
{\tt http://www-fourier.ujf-grenoble.fr/\~{}demailly/books.html}&

\Bibitem [Dem21]&Demailly, J.-P.:& hermitian-Yang-Mills approach to
the conjecture of Griffiths on the positivity of ample vector bundles,&
Mat.\ Sb.\ {\bf 212}, Number 3 (2021) 39--53&

\Bibitem [DPS94]&Demailly, J.-P., Peternell, Th., Schneider, M.:& Compact
complex manifolds with numerically effective tangent bundles,&
J.~Algebraic Geometry {\bf 3} (1994) 295--345&

\Bibitem [Don85]&Donaldson, S.:& Anti-self-dual Yang-Mills connections over
complex algebraic surfaces and stable vector bundles,& Proc.\ London Math.\
Soc. {\bf 50} (1985), 1--26&

\Bibitem [Gri69]&Griffiths, P.A:& hermitian differential geometry,
Chern classes and positive vector bundles,& Global Analysis, papers in
honor of K. Kodaira, Princeton Univ. Press, Princeton (1969),
181--251&

\Bibitem [Kod54]&Kodaira, K.:& On K\"ahler varieties of restricted type,&
Ann.\ of Math.\ {\bf 60} (1954) 28--48&

\Bibitem [LSY13]&Liu K., Sun X., Yang, X.:& Positivity and vanishing theorems
for ample vector bundles,& J.~Alg.\ Geom.\ {\bf 22} (2013) 303–-331&

\Bibitem [MoT07]&Mourougane, C., Takayama, S.:& Hodge metrics and
positivity of direct images,& J.~reine angew.\ Math.\ {\bf 606}
(2007), 167--179&

\Bibitem [Nak55]&Nakano, S.:& On complex analytic vector bundles,&
J.~Math.\ Soc.\ Japan {\bf 7} (1955) 1--12&

\Bibitem [Pin20]&Pingali, V.P.:& A vector bundle version of the Monge-Amp\`ere
equation,& Adv.\ in Math.\ {\bf 360} (2020), 40 pages,
https://doi.org/10.1016/j.aim.2019.106921&

\Bibitem [Pin21]&Pingali, V.P.:& A note on Demailly's approach towards a
conjecture of Griffiths,& C.~R.\ Math.\ Acad.\ Sci.\ Paris, {\bf 359}
(2021) 501--503&

\Bibitem [UhY86]&Uhlenbeck, K., Yau, S.T.:& On the existence of 
Hermitian-Yang-Mills connections in stable vector bundles,& Comm.\
Pure and Appl.\ Math.\ {\bf 39} (1986) 258--293&

\Bibitem [Yau78]&Yau, S.T.:& On the Ricci curvature of a complex K\"ahler
manifold and the complex Monge-Amp\`ere equation I,& Comm.\ Pure and Appl.\
Math.\ {\bf 31} (1978) 339--411&
\vskip15pt

\parindent=0cm
(Version of January 9, 2022, printed on \today, \timeofday)
\bigskip

\noindent
Jean-Pierre Demailly\\
Universit\'e Grenoble Alpes, Institut Fourier (Math\'ematiques)\\
UMR 5582 du C.N.R.S., 100 rue des Maths, 38610 Gi\`eres, France\\
{\it e-mail:}\/ jean-pierre.demailly@univ-grenoble-alpes.fr

\end{document}